\documentclass[smallextended,envcountsect]{svjour3}
\spnewtheorem{Assumption}{Assumption}[section]{\bfseries}{\itshape}
\smartqed
\usepackage{graphicx,amsmath,amsfonts,color}
\usepackage{setspace}
\usepackage{pdfsync}
\journalname{JOTA}

\def\ben{\begin{eqnarray}}
\def\een{\end{eqnarray}}
\def\be{\begin{eqnarray}}
\def\ee{\end{eqnarray}}

\def\beno{\begin{eqnarray*}}
\def\eeno{\end{eqnarray*}}

\DeclareMathOperator*{\esssup}{ess\,sup}
\allowdisplaybreaks[1]

\begin{document}

\title{On Dynamic Programming Principle for Stochastic Control under Expectation Constraints}

%\subtitle{Using  the  LaTex Template}

\titlerunning{J Optim Theory Appl} % if too long for running head

\author{Yuk-Loong  Chow,\ \ Xiang Yu,\ \ Chao Zhou}

\institute{Communicated by Lars Gr\"{u}ne\\
\ \\
\ \\
Yuk-Loong Chow \at
             School of Mathematics, Sun Yat-Sen University \\
              Guangzhou, China\\
              zhouyulong@mail.sysu.edu.cn
           \and
           Xiang Yu, Corresponding author \at
              Department of Applied Mathematics, The Hong Kong Polytechnic University \\
              Hung Hom, Kowloon, Hong Kong\\
              xiang.yu@polyu.edu.hk
           \and
           Chao Zhou \at
              Department of Mathematics, National University of Singapore \\
              Singapore\\
              matzc@nus.edu.sg             
}

\date{  }

\maketitle

\begin{abstract}
This paper studies the dynamic programming principle using the measurable selection method for stochastic control of continuous processes. The novelty of this work is to incorporate intermediate expectation constraints on the canonical space at each time $t$. Motivated by some financial applications, we show that several types of dynamic trading constraints can be reformulated into expectation constraints on paths of controlled state processes. Our results can therefore be employed to recover the dynamic programming principle for these optimal investment problems under dynamic constraints, possibly path-dependent, in a non-Markovian framework.
\end{abstract}

%All acknowledgements should be placed in the back of the paper after Conclusions..

\section{Introduction}\label{sec-int}
Stochastic control problems have been actively studied and widely applied in many areas since the 1970s, especially in quantitative finance. To determine the optimal control process and compute the corresponding value function, the celebrated dynamic programming principle (DPP) plays a crucial role. It provides a convenient way to handle a global optimization problem by solving a series of recursive local optimization problems, which is closely related to the tower property of conditional expectations. In the Markovian framework, it has been shown in the vast literature that DPP allows one to relate the value function of a control problem to a solution of the associated Hamilton-Jacob-Bellman (HJB) equation in a classical sense or a viscosity sense. Numerical approximations can also come into play in various models. See \cite{ET13a,ET13b} for the comprehensive review on the topic of stochastic control.

The growing complexity of financial markets motivates the continuing development of stochastic control theory. In particular,  stochastic control problems subject to various constraints have been extensively investigated in different contexts. The existing literature is far-reaching, see for instance \cite{Bouchard-Elie-Imbert,Bouchard-Elie-Touzi,Bouchard-Nutz,Soner1,Soner2,Soner-Touzi,Soner-Touzi2} that discuss DPP for stochastic control with constraints and also \cite{BayChrisMil,BayCoxSto,BayYao,SKad} that address DPP issues for optimal stopping with constraints.
To avoid measurability issues of the value function and technical challenges in measurable selection arguments,
one can introduce some measurable test functions based on the value function and develop the weak version of DPP instead.
Note that a weak DPP is sufficient to verify that the value function corresponds to a viscosity solution to some HJB equation, see \cite{Bouchard-Touzi}. In a similar fashion, \cite{Bouchard-Nutz} further establishes a weak DPP for control problems under the expectation constraint at terminal time in a Markovian setting.

Dynamic trading constraints arise naturally in many financial applications. In particular, expectation constraints at each intermediate time create new challenges in the control problem, and it is still an open problem whether the strong version of DPP holds or not, especially in the non-Markovian setting. This paper aims to fill this gap and proves the strong DPP by applying some measurable selection arguments. For a given control process, the resulting state process induces a new probability measure. Consequently, to find an optimal control is equivalent to optimize the probability measures induced by all admissible controls. Along this line, \cite{Nutz-Handel} constructed a sublinear expectation which satisfies the tower property also known as time-consistency. Regarding probabilities therein as controls and the sublinear expectation as the value function, the tower property is essentially equivalent to DPP. To handle the expectation constraints at terminal time, \cite{Bouchard-Nutz} introduces some auxiliary martingales to aid the proof of DPP. We generalize this idea further by considering some auxiliary supermartingales and establish the desired DPP in the model when intermediate expectation constraints are binding at all time. Moreover, our general setting allows the admissible probabilities and the value function to be path-dependent so that DPP in this paper can be applied in some path-dependent stochastic control problems.

Another contribution of this paper is to show that several dynamic trading constraints such as (but not limited to) the following types can be transformed into the intermediate expectation constraint:\\
\indent$\mathrm{(i)}$ State constraint: it requires the controlled wealth process to stay in a prescribed region, see some examples in \cite{Soner1,Soner2,Ishii-Koike,Lasry-Lions,Katsoulakis}.\\
\indent$\mathrm{(ii)}$ Floor constraint: it requires the controlled wealth process to stay above a benchmark stochastic process, see \cite{EIKaroui-Jeanblanc-Lacoste}.\\
\indent$\mathrm{(iii)}$ Drawdown constraint: it requires the controlled wealth process to stay above a fraction of its own running maximum process, see \cite{Grossman-Zhou} and \cite{Elie-Touzi}.\\
\indent$\mathrm{(iv)}$ Target constraint: it studies a reachability set, namely, the collection of all initial data such that the state process can be driven into a target set at a specified time. See some examples in \cite{Soner-Touzi} and \cite{Bouchard-Elie-Imbert,Bouchard-Elie-Touzi,Soner-Touzi2}.\\
\indent$\mathrm{(v)}$ Quantile hedging: it mandates the wealth process to stay in a given set with the probability greater than some specified level, see \cite{Follmer-Leukert} and \cite{Follmer-Leukert2}.\\
Our main results can also recover DPP in other control applications with constraints beyond quantitative finance such as aircraft abort landing in \cite{PVinter-1} and \cite{PVinter} and energy and resources management in \cite{Kharroubi}. See also \cite{ST-1} and \cite{BPZ16} on general stochastic control problems with similar constraints.

The rest of the paper is organized as follows. Sect. \ref{sec-main} introduces some mathematical preliminaries. An abstract version of DPP is first developed therein. Sect. \ref{sec-cont} establishes the main result of this paper, namely, DPP for stochastic control problems under intermediate expectation constraints. Sect. \ref{sec-app} presents some applications in optimal investment and hedging problems under various dynamic trading constraints. At last, we give the conclusion in Sect. \ref{section_conc}.

\section{General Framework}\label{sec-main}
\subsection{\textbf{Notations and Preliminaries}}
Let $\Omega =  \{\omega \in C([0,T];\mathbb{R}^{d}):\omega_{0} = 0\}$ be the canonical space of continuous paths equipped with the topology of uniform convergence under the norm $\|\omega\|_{\infty}:=\sup_{t\in[0,T]}\|\omega_t\|$. $P_{0}$ is the Wiener measure on $\Omega$ and $B$ is the canonical process $B_{t}(\omega)=\omega_{t}$.

Let $\mathbb{F} = (\mathcal{F}_{t})_{t \in [0,T]}$ be the canonical filtration generated by $B$ and let $\mathbb{F}_+=(\mathcal{F}_{t}^+)_{t \in [0,T]}$ be the right limit of $\mathbb{F}$ with $\mathcal{F}_{t}^+:=\cap_{s>t}\mathcal{F}_{s}$ for all $t \in [0,T)$ and $\mathcal{F}_{T}^+=\mathcal{F}_{T}$. Furthermore, $\mathfrak{P}(\Omega)$ denotes the set of all probability measures on $(\Omega,\mathcal{F})$, equipped with the topology of weak convergence. Given two paths $\bar{\omega}$ and $\omega$, their concatenation at $t$ is defined by
$$ (\bar{\omega}\otimes_t\omega)_s:=\bar{\omega}_s\mathbf{1}_{[0,t[}(s)+(\bar{\omega}_t+\omega_s-\omega_t)\mathbf{1}_{[t,T]}(s),\ s\in[0,T].
$$
Let $E^{P}[\cdot]$ denote the expectation under probability $P$, and we simply write $E[\cdot]$ if $P=P_{0}$. We define integrals of measurable functions $\xi$ with values in the extended real line $\bar{\mathbb{R}}=[-\infty,\infty]$ and set $E^{P}[\xi]:=E^{P}[\xi^{+}]-E^{P}[\xi^{-}]$
if $E^{P}[\xi^{+}]$ or $E^{P}[\xi^{-}]$ is finite, and define $E^{P}[\xi]:= -\infty ~ \text{if} ~E^{P}[\xi^{+}]=E^{P}[\xi^{-}]=+\infty$. In accordance with the convention, we adopt $\sup \emptyset = -\infty$ (resp. $\inf \emptyset = \infty$). Let $\mathcal{T}$ be the collection of all $\mathbb{F}-$stopping times taking values in $[0,T]$.

As in Chapter 1 of \cite{Soner2} and Chapter 7 of \cite{Bertsekas-Shreve}, for any $P \in \mathfrak{P}(\Omega), \tau \in \mathcal{T},$ there exists a regular conditional probability distribution $\{P^{\omega}_{\tau}\}_{\omega\in\Omega}$ of $P$ given $\mathcal{F}_{\tau}$. That is, $P^{\omega}_{\tau}\in \mathfrak{P}(\Omega)$ for each $\omega$, and $\omega \in \Omega \rightarrow P^{\omega}_{\tau}(A) \in [0,1]$ is $\mathcal{F}_{\tau}$-measurable for any $A\in\mathcal{F}$ and
$E^{P^{\omega}_{\tau}}[\xi]=E^{P}[\xi|\mathcal{F}_{\tau}](\omega) ~\text{for} ~ P\text{-a.e.}~\omega \in \Omega,$
whenever $\xi$ is $\mathcal{F}$-measurable and bounded. Moreover, $P^{\omega}_{\tau}$ can be chosen to be concentrated on the set of paths that coincide with $\omega$ up to time $\tau(\omega).$ 
Given $P \in \mathfrak{P}(\Omega)$ and a family  $(Q^{\omega})_{\omega\in\Omega}$ such that $\omega \in \Omega \rightarrow Q^{\omega} \in \mathfrak{P}(\Omega)$ is $\mathcal{F}_{\tau}$-measurable with $Q^{\omega}(\Omega^{\omega}_{\tau})=1$ for all $\omega \in \Omega$, one can define a concatenated probability measure $P\otimes_{\tau}Q^{\cdot}$ by
$P\otimes_{\tau}Q^{\cdot}(A):=\int_{\Omega}Q^{\omega}(A)P(d\omega), ~~ \forall A \in \mathcal{F}.$ As in \cite{Bertsekas-Shreve}, a subset of a Borel space is called analytic if it is the image of a Borel subset of another Borel space under a Borel-measurable function. 

\subsection{\textbf{Problem Formulation and Main Results}}
For each $(t,\omega)\in[0,T]\times\Omega,$ we consider an adapted set $\mathcal{P}(t,\omega) \subseteq \mathfrak{P}(\Omega)$ such that $\mathcal{P}(t,\omega)  = \mathcal{P}(t,\tilde{\omega})$, if $\omega = \tilde{\omega}$ on $[0,t]$ and it is assumed that $\mathcal{P}(t,\omega) \neq \emptyset$. Let $\Omega^{\omega}_{t}:=\{\omega^{\prime}\in\Omega: \omega^{\prime}=\omega ~\text{on} ~[0, t]\},$ then $P(\Omega^{\omega}_{t})=1,$ for any
$P \in \mathcal{P}(t,\omega)$. If $\tau$ is a stopping time, we denote $\mathcal{P}(\tau,\omega) = \mathcal{P}(\tau(\omega),\omega)$. For $t=0$, we note that $\mathcal{P}_0:=\mathcal{P}(0,\omega)$ for any $\omega\in\Omega$.

\begin{Assumption} \label{ass:0}
We assume that every $P\in\mathcal{P}_0$ satisfies the Blumenthal zero-one law
and the martingale representation property.
\end{Assumption}

\begin{remark}
By Remark 2.4 in \cite{STZ}, it follows from Blumenthal zero-one law that $E^P[\xi\mid\mathcal F_t]=E^P[\xi\mid\mathcal F_t^+]$, $P$-a.s. for any $t\in[0,1]$ and $P$-integrable random variable $\xi$. That is, any $\mathcal F_t^+$-measurable random variable has an $\mathcal F_t$-measurable $P$-modification. One typical example that satisfies Assumption \ref{ass:0} is the space of probability measures induced by strong solutions of controlled SDEs with invertible volatility coefficients. Our setting therefore includes some natural examples that stem from financial applications.
\end{remark}

Throughout the paper, we fix $\xi: \Omega \rightarrow \bar{\mathbb{R}}$ and $\eta_t: \Omega \rightarrow \bar{\mathbb{R}}$ for each $t\in[0,T]$.
In our control problem, $\xi$ is used to model the reward and $\eta_t$ is used to model the intermediate constraint.
It is assumed that $\xi$ is upper semi-analytic and $\eta_t$ is lower semi-analytic for each $t\in[0,T]$ and the process $(\eta_t)_{t\in[0,T]}$ has lower semi-continuous paths. For a given constraint level $m$, we define the set
\ben \label{expectation-constraints}
\mathcal{P}(t,\omega,m) :=\{P\in \mathcal{P}(t,\omega): E^P[\eta_s]\leq m,\ \forall t\leq s\leq T\}.
\een
Let us consider the control problem with intermediate expectation constraints that
\begin{equation}\label{value-function}
V(t,\omega,m) = \sup_{P \in \mathcal{P}(t,\omega,m)} E^{P}[\xi].
\end{equation}
Note that the formulation \eqref{value-function} is general enough to cover many control problems including running rewards and exponential discounting.
The present paper aims to verify DPP for problem \eqref{value-function}, which relates the value $V(t)$ to the conditional expectation of the later value
$V(\tau)$ at some stopping time $\tau \geq t$. At a later time $\tau$, the state $\omega$ incurs more realization and the constraint level may also change. We therefore need to randomize the constraints by introducing some auxiliary supermartingales, which are motivated by the auxiliary martingales used in \cite{Bouchard-Nutz}.

\begin{definition}\label{defM}
For each $P \in \mathcal{P}(t,\omega, m)$, let $\mathcal{M}^{+}_{t,\omega,m}(P)$ be the collection of all supermartingales on $[t, T] \times \Omega \rightarrow \bar{\mathbb{R}}$ such that\\
$\mathrm{(i)}$ $M_{t} \leq m$, $P$-a.s.;\\
$\mathrm{(ii)}$ $E^{P} [M_{s_{2}}|\mathcal{F}_{s_{1}}] \leq M_{s_{1}}$, $P$-a.s., for any $s_{1}, s_{2} \in [t,T]$ with $ s_{1} \leq s_{2}$;\\
$\mathrm{(iii)}$ $M_{s} \geq \eta_s$, $P$-a.s. for any $t\leq s\leq T$.\\
\end{definition}

First, we need to verify that the set $\mathcal{M}^{+}_{t,\omega,m}(P)$ is not empty. One can consider the set $\mathbf{V}_{t, m}$ of martingales starting from $m$ at time $t$ that 
\begin{align*}
\mathbf{V}_{t, m}:= \{M^{t,m,\alpha}: M^{t,m,\alpha}= m+\int_t^{\cdot}\alpha_s dW^P_s\},
\end{align*} 
for some $\alpha \in \textbf{H}_2$ from the martingale representation of $(E^P[\eta_s | \mathcal{F}_u])_{u\in[t,s]}$ \text{under} $P$, where $W^P$ is a $P$-Brownian motion and $\textbf{H}_2$ denotes the set of all adapted and square integrable processes. However, for {\it a priori} $s\geq t$, we can only get the existence of $\alpha^{(s)} \in \textbf{H}_2$ such that the controlled martingale $M_{\cdot}^{t,m,\alpha^{(s)}}$ satisfies Definition \ref{defM} up to time $s$. More efforts are required in order to ensure the existence of a controlled martingale independent on the time $s$ using the idea of aggregation. Let $\mathcal T^{0}$ denote the set of $\mathbb F$-stopping times $\tau$ such that $\tau \in [0,T]$ a.s. For $\theta$ in $\mathcal T^{0}$, $\mathcal T^{\theta}$  is the set of stopping times
$\tau\in \mathcal T^{0}$ such that $\theta \leq \tau \leq T$, $P$-a.s.. The next result confirms the existence of such an aggregated auxiliary supermartingale.

\begin{lemma}\label{keyres}
Under Assumption \ref{ass:0}, the set $\mathcal{M}^{+}_{t,\omega,m}(P)$ is non-empty.
\end{lemma}
%\noindent{\it Proof}
%For each $\sigma \in \mathcal{T}^0,$ we define the $\mathcal{F}_\sigma$-measurable random variable
%\begin{equation}\label{value}
%V_{\sigma}:=  \esssup_{\tau \in \mathcal{T}^{\sigma}} E^P\left[\eta_\tau| \mathcal{F}_\sigma\right].
%\end{equation}
%By Theorem VI-48 in \cite{Dell-1} {\color{blue}for continuous processes}, the family of supermartingales $\{V_{\sigma}:\sigma \in \mathcal{T}^0\}$ can be aggregated by an optional process $(V_t)$ admitting the Mertens decomposition that {\color{blue}$V_t:=N_t-A_t$}, where $N$ is a continuous square integrable martingale and $A$ is a non-decreasing continuous predictable process such that $A_0=0$. As for all $s\in[t,T]$, it holds that ${E}^P\left[\eta_s| \mathcal{F}_t\right] \leq m $ $P$-a.s., we get that $\underset{s \in [t,T]}{\esssup}\,E^P\left[\eta_s| \mathcal{F}_t\right] \leq m$. In view of the definition of $V$ (see \eqref{value}), we obtain
%\begin{equation}\label{tau}
%V_t=N_t-A_t\leq m {\rm\,\, P-a.s.}
%\end{equation}
%For fixed $s \geq t$, we have $\eta_s=E^P\left[\eta_s| \mathcal{F}_s\right] \leq\underset{u \in [s,T]}{\esssup}\, E^P\left[\eta_u| \mathcal{F}_s \right]=V_s {\rm \,\, P-a.s.}$
%This observation together with \eqref{tau} lead to
%$$\eta_s \leq N_s-A_s=N_t-A_t+\int_t^s \alpha_udW^P_u-A_s+A_t.$$
%Using the above inequality, \eqref{tau} and the fact that the processes $A$ is non-decreasing, we obtain that $\eta_s \leq M_s^{t,m,\alpha} {\rm \,\, P-a.s.}$ and our claim holds.
%\qed
{\it Proof}
For each $\sigma \in \mathcal{T}^0,$ we define the $\mathcal{F}_\sigma$-measurable random variable
\begin{equation}\label{value}
%V_{\sigma}:=  \esssup_{\tau \in \mathcal{T}^{\sigma}} E^P\left[\eta_\tau| \mathcal{F}_\sigma\right].
\mathcal V_{\sigma}:=  \esssup_{\tau \in \mathcal{T}^{\sigma}} E^P\left[\eta_\tau| \mathcal{F}_\sigma\right]\ \ \ P-a.s.
\end{equation}
By classical results of the general theory of processes (see e.g. \cite{Dell-1}) and Remark 2.1, %By Theorem VI-48 in \cite{Dell-1},
the family of supermartingales $\{V_{\sigma}:\sigma \in \mathcal{T}^0\}$ can be aggregated by an optional process $(\mathcal V_t)$ admitting the Mertens decomposition that 
\begin{align*}
\mathcal V_t:=N_t-A_t-C_{t^-},
\end{align*}
where $N$ is a square integrable martingale, $A$ is a non-decreasing RCLL predictable process such that $A_0=0$ and $C$ is a non-decreasing right-continuous adapted process, purely discontinuous satisfying $C_{0^-}=0$. %\\

As for all $s\in[t,T]$, it holds that ${E}^P\left[\eta_s| \mathcal{F}_t\right] \leq m $, $P$-a.s., which further entails that $\underset{s \in [t,T]}{\esssup}\,E^P\left[\eta_s| \mathcal{F}_t\right] \leq m$. In view of the definition of $\mathcal V$ (see \eqref{value}), we obtain
\begin{equation}\label{tau}
\mathcal V_t=N_t-A_t-C_{t^-} \leq m, \,\, P{\rm -a.s.}
\end{equation}
For fixed $s \geq t$, we have $\eta_s=E^P\left[\eta_s| \mathcal{F}_s\right] \leq\underset{u \in [s,T]}{\esssup}\, E^P\left[\eta_u| \mathcal{F}_s \right]=\mathcal V_s, \,\, P{\rm -a.s.}$
This observation, together with \eqref{tau}, leads to
$$\eta_s \leq N_s-A_s-C_{s^-}=N_t-A_t-C_{t^-}+\int_t^s \alpha_udW^P_u-A_s+A_t-C_{s^-}+C_{t^-}.$$
Using the above inequality, \eqref{tau} and the fact that the processes $A$ and $C$ are non-decreasing, we obtain $\eta_s \leq M_s^{{t,m,\alpha}}$, $P$-a.s. and our claim holds.
%Therefore the claims holds that the set $\mathcal{M}^{+}_{t,\omega,m}(P)$ is non-empty.
\qed

Given $\tau \in \mathcal{T}^{t}$, $\omega \in \Omega$, $m\in \mathbb{R}$, $P\in \mathcal{P}(t,\omega,m)$, $M \in \mathcal{M}^{+}_{t,\omega,m}(P)$, we set $\mathcal{P}(\tau,\omega,M_{\tau}) := \mathcal{P}(\tau(\omega),\omega,M_{\tau}(\omega)), V(\tau,\omega,M_{\tau}) := V(\tau(\omega),\omega, M_{\tau}(\omega))$.
The following conditions are required in our main results.
\begin{Assumption} \label{ass:1}
Let $(t,\bar{\omega})\in[0,T]\times\Omega$ and $\tau\in\mathcal{T}^{t}$, for any $P\in\mathcal{P}(t,\bar{\omega})$, we assume:\\
$\mathrm{(i)}$ Measurability: the graph $[[\mathcal{P}]]:= \{(t,\omega,Q): (t,\omega)\in[0,T]\times\Omega, Q \in\mathcal{P}(t, \omega)\}$ is an analytic subset of $[0,T]\times\Omega\times\mathfrak{P}(\Omega)$.\\
$\mathrm{(ii)}$ Invariance: there is a family of regular conditional probability distribution $(P^{\omega}_{\tau})_{\omega\in\Omega}$ of $P$ given $\mathcal{F}_{\tau}$ such that $P^{\omega}_{\tau} \in \mathcal{P}(\tau,\omega)$ for $P$-a.e. $\omega \in \Omega$.\\
$\mathrm{(iii)}$ Stability under pasting: let $(Q^{\omega})_{\omega\in\Omega}$ be such that $\omega \rightarrow Q^{\omega}$ is $\mathcal{F}_{\tau}-$measurable and $Q^{\omega}\in\mathcal{P}(\tau,\omega)$ for $P$-a.e. $\omega \in \Omega$, then $P\otimes_{\tau} Q^{\cdot}\in \mathcal{P}(t,\bar{\omega}).$
\end{Assumption}

\begin{remark} Assumption \ref{ass:1} is essentially the same as Assumption 2.1 in \cite{Katsoulakis}, in which $G-$expectations and random $G-$expectations have been considered as typical examples.

\end{remark}

Our first main result gives DPP in a measurable selection setting.
\begin{theorem} \label{thm:main}
Under Assumption \ref{ass:0} and Assumption \ref{ass:1}, the value function $V$ defined by \eqref{value-function} satisfies
\begin{equation} \label{eq2.1}
\begin{split}
V(t,\bar{\omega},m) & = \sup_{P\in\mathcal{P}(t,\bar{\omega},m)}\sup_{M\in\mathcal{M}^{+}_{t,\bar{\omega},m}(P)}E^{P}[V(\tau,\bar{\omega}\otimes_t\omega,M_{\tau})] \\
 & = \sup_{P\in\mathcal{P}(t,\bar{\omega},m)}\inf_{M\in\mathcal{M}^{+}_{t,\bar{\omega},m}(P)}E^{P}[V(\tau,\bar{\omega}\otimes_t\omega,M_{\tau})].
\end{split}
\end{equation}
\end{theorem}

\begin{remark} The DPP result \eqref{eq2.1} involves both supremum and infimum over $\mathcal{M}^{+}_{t,\bar{\omega},m}(P)$, which manifests the auxiliary role of
$\mathcal{M}^{+}_{t,\bar{\omega},m}(P)$. In other words, any element in $\mathcal{M}^{+}_{t,\bar{\omega},m}(P)$ plays the same role.
Intuitively,  under the optimal probability $P_{*}$, only the expectation matters and different choices of super-martingales in $\mathcal{M}^{+}_{t,\bar{\omega},m}(P_{*})$ will lead to the same expectation.
\end{remark}

The proof of Theorem \ref{thm:main} requires several auxiliary results. First, by \cite{Bertsekas-Shreve}, Corollary 7.48.1, we have
\begin{lemma} \label{lem:1}
If $\xi$ is upper semi-analytic, the function $P\in\mathfrak{P}(\Omega)\rightarrow E^{P}[\xi]\in\bar{\mathbb{R}}$ is upper semi-analytic. Similarly, for each $t\leq s\leq T$, if $\eta_s$ is lower semi-analytic, the function $P\in\mathfrak{P}(\Omega)\rightarrow E^{P}[\eta_s]\in\bar{\mathbb{R}}$ is lower semi-analytic.
\end{lemma}

Similar to the proof of Lemma $4.1$ and Proposition $4.1$ in \cite{BayMih}, we have the following result.
\begin{lemma}\label{lema:1-1}
As $(\eta_s)_{s\in[t,T]}$ has lower semi-continuous paths, there exists a countable subset $\mathcal{S}:=\{s_n: t\leq s_n\leq T\}_{n\in\mathbb{N}}$ such that the supremum can be achieved, i.e.,
$
\sup_{s\in[t,T]}E^{P}[\eta_s]=\sup_{s_n\in\mathcal{S}}E^{P}[\eta_{s_n}].
$
\end{lemma}

Next, we give a measurability result.
\begin{lemma} \label{lem:2} Let $D := \{(t,\omega,m,P): (t,\omega,m)\in[0,T]\times\Omega\times\mathbb{R}, P \in\mathcal{P}(t, \omega,m)\}$, then $D$ is an analytic subset of  $[0,T]\times\Omega\times\mathbb{R}\times\mathfrak{P}(\Omega)$.
\end{lemma}
{\it Proof}
 Observe that
\begin{equation}
\begin{split}
D = &\left\{(t,\omega,m,P): (t,\omega,m)\in[0,T]\times\Omega\times\mathbb{R}, P \in\mathcal{P}(t, \omega)\right\} \cap\\
 &\left\{(t,\omega,m,P): L(s,m,P) \leq 0,\ \forall t\leq s\leq T\right\}\\
=&\left\{(t,\omega,m,P): (t,\omega,m)\in[0,T]\times\Omega\times\mathbb{R}, P \in\mathcal{P}(t, \omega)\right\} \cap\\
& \{(t,\omega,m,P): \sup_{s\in[0,T]}L(s,m,P) \leq 0\}.\nonumber
\end{split}
\end{equation}
The first term is analytic by item $\mathrm{(i)}$ of Assumption \ref{ass:1}. 

We claim that the second term is also analytic. Firstly, Lemma \ref{lem:1} gives that $L(s,m,P)$ is lower semi-analytic for each fixed $t\leq s\leq T$ and Lemma \ref{lema:1-1} asserts that $\sup_{s\in[0,T]}L(s,m,P)=\sup_{s_n\in\mathcal{S}}L(s_n,m,P)$. Secondly, Lemma $7.30$ of \cite{Bertsekas-Shreve} entails that $\sup_{s\in[0,T]}L(s,m,P)$ is also lower semi-analytic because $\sup_{s_n\in\mathcal{S}}L(s_n,m,P)$ is lower semi-analytic. The claim therefore holds and we get that $D$ is analytic.
\qed

To simplify the notation, let us set $X = [0,T]\times\Omega\times\mathbb{R},$ and define $\text{proj}_{X}(D) = \{(t,\omega,m): (t,\omega,m,P) \in D\}$. It then follows that 
\begin{align*}
\text{proj}_{X}(D) = \{(t,\omega,m): \mathcal{P}(t,\omega,m) \neq \emptyset\}.
\end{align*}

\begin{lemma} \label{lem:3} The value function $V:\text{proj}_{X}(D)\rightarrow\bar{\mathbb{R}}$ is  upper semi-analytic. Moreover, for every $\epsilon > 0$, there exists an analytically measurable function $\varphi_{\epsilon}: \text{proj}_{X}(D)\rightarrow\mathfrak{P}(\Omega)$ such that for every $(t,\omega,m)\in \text{proj}_{X}(D),$ one has $(t,\omega,m, \varphi_{\epsilon}(t,\omega,m)) \in D$ and
\begin{equation}
  E^{\varphi_{\epsilon}(t,\omega,m)}[\xi]\geq \begin{cases}
               V(t,\omega,m)-\epsilon, &  V(t,\omega,m)<\infty,\\
               \epsilon^{-1},             &  V(t,\omega,m)=\infty.
            \end{cases}
\end{equation}
\end{lemma}
{\it Proof}
Thanks to Lemma \ref{lem:1} and Lemma \ref{lem:2}, it is easy to conclude Lemma \ref{lem:3} by Proposition 7.47 on page 179 and Proposition 7.50 on page 178 of \cite{Bertsekas-Shreve}.
\qed

\noindent{{\it Proof of Theorem \ref{thm:main}}}:
\textit{Step 1:} Let us first show one direction of \eqref{eq2.1} that
\begin{equation} \label{eq:2.2}
V(t,\bar{\omega},m)  \leq \sup_{P\in\mathcal{P}(t,\bar{\omega},m)}\inf_{M\in\mathcal{M}^{+}_{t,\bar{\omega},m}(P)}E^{P}[V(\tau,\bar{\omega}\otimes_t\omega,M_{\tau})].
\end{equation}
Fix $P\in\mathcal{P}(t,\bar{\omega},m)$ and $M\in\mathcal{M}^{+}_{t,\bar{\omega},m}(P)$.  By Assumption \ref{ass:1} $\mathrm{(ii)}$, there exists a family of regular conditional probability distributions $(P^{\bar{\omega}\otimes_t\omega}_{\tau})$ of $P$ given $\mathcal{F}_{\tau}$ such that $P^{\bar{\omega}\otimes_t\omega}_{\tau} \in \mathcal{P}(\tau,\bar{\omega}\otimes_t\omega)$ for $P$-a.e. $\omega\in \Omega$. We then claim that $P^{\bar{\omega}\otimes_t\omega}_{\tau} \in \mathcal{P}(\tau,\bar{\omega}\otimes_t\omega,M_{\tau})$ for $P$-a.e. $\omega \in \Omega$. To see this, for $P$-a.e. $\omega\in\Omega$ and $\tau\leq \rho \leq T$, we have
\begin{equation} \label{eq:2.3}
E^{P^{\bar{\omega}\otimes_t\omega}_{\tau}}[\eta_{\rho}]
  = E^{P}[\eta_{\rho}|\mathcal{F}_{\tau}](\bar{\omega}\otimes_t\omega)
  \leq E^{P}[M_{\rho}|\mathcal{F}_{\tau}](\bar{\omega}\otimes_t\omega)
  \leq M_{\tau}(\bar{\omega}\otimes_t\omega).
\end{equation}
The claim therefore holds. It follows that for $P$-a.e. $\omega\in\Omega$, one has
\begin{equation} \label{eq:2.4}
E^{P}[\xi|\mathcal{F}_{\tau}](\bar{\omega}\otimes_t\omega)
  = E^{P^{\bar{\omega}\otimes_t\omega}_{\tau}}[\xi]
  \leq V(\tau,\bar{\omega}\otimes_t\omega,M_{\tau}).
\end{equation}
Taking $P(d\omega)$-expectations, we obtain
$E^{P}[\xi] \leq E^{P}[V(\tau,\bar{\omega}\otimes\omega,M_{\tau})],$ which gives
$E^{P}[\xi] \leq \inf_{M\in\mathcal{M}^{+}_{t,\bar{\omega},m}(P)}E^{P}[V(\tau,\bar{\omega}\otimes\omega,M_{\tau})].$
The inequality (\ref{eq:2.2}) follows by taking the supremum over $\mathcal{P}(t,\bar{\omega},m).$\\
\noindent
\textit{Step 2:} We then turn to prove the opposite direction that
\begin{equation} \label{eq:2.5}
V(t,\bar{\omega},m)  \geq \sup_{P\in\mathcal{P}(t,\bar{\omega},m)}\sup_{M\in\mathcal{M}^{+}_{t,\bar{\omega},m}(P)}E^{P}[V(\tau,\bar{\omega}\otimes_{t}\omega,M_{\tau})].
\end{equation}
Fix $\epsilon > 0, P\in\mathcal{P}(t,\bar{\omega},m)$, and take an arbitrary $M\in\mathcal{M}^{+}_{t,\bar{\omega},m}(P)$. As the composition of universally measurable functions is universally measurable, the map $\omega \in \Omega \rightarrow \varphi_{\epsilon}(\tau(\bar{\omega}\otimes_t\omega),\bar{\omega}\otimes_t\omega,M_{\tau}(\bar{\omega}\otimes_t\omega)) \in \mathfrak{P}(\Omega)$ is $\mathcal{F}^{*}_{\tau}$-measurable by the universally measurable extension of Galmarino's test, see lemma 2.5 in \cite{Nutz-Handel}. Therefore, there exists an $\mathcal{F}_{\tau}$-measurable kernel $Q_{\epsilon}: \Omega \rightarrow \mathfrak{P}(\Omega)$ such that $Q_{\epsilon}^{\omega} = \varphi_{\epsilon}(\tau(\bar{\omega}\otimes_t\omega),\bar{\omega}\otimes_t\omega,M_{\tau}(\bar{\omega}\otimes_t\omega))$
for $P$-a.e. $\omega\in\Omega$. Again by Assumption \ref{ass:1}(2) and equation (\ref{eq:2.3}), we have $\mathcal{P}(\tau,\omega,M_{\tau}) \neq \emptyset$ for $P$-a.e. $\omega\in\Omega$. Thus by Lemma \ref{lem:3}, for $P$-a.e. $\omega\in\Omega$, we have $Q_{\epsilon}^{\omega}\in \mathcal{P}(\tau,\bar{\omega}\otimes_t\omega,M_{\tau})$ and
\[
  E^{Q^{\omega}_{\epsilon}}[\xi] \geq \begin{cases}
               V(\tau,\bar{\omega}\otimes_t\omega,M_{\tau})-\epsilon, &  V(\tau,\bar{\omega}\otimes_t\omega,M_{\tau})<\infty,\\
               \epsilon^{-1},             &  V(\tau,\bar{\omega}\otimes_t\omega,M_{\tau})=\infty.
            \end{cases}
\]
It yields that $P\otimes_{\tau}Q_{\epsilon}^{\cdot} \in \mathcal{P}(t,\bar{\omega})$ by item $\mathrm{(iii)}$ of Assumption \ref{ass:1}. We now claim that $P\otimes_{\tau}Q_{\epsilon}^{\cdot} \in \mathcal{P}(t,\bar{\omega},m)$. To see this, for any $t \leq \rho \leq T$, we have 
\beno
E^{P\otimes_{\tau}Q_{\epsilon}^{\cdot}}[\eta_{\rho}]
 &=& E^{P}[E^{Q_{\epsilon}^{\cdot}}[\eta_{\rho}\textbf{1}_{\rho > \tau}]] + E^{P}[E^{Q_{\epsilon}^{\cdot}}[\eta_{\rho}\textbf{1}_{\rho \leq \tau}]]
  \\&\leq& E^{P}[M_{\tau}\textbf{1}_{\rho > \tau}]
+E^{P}[M_{\rho}\textbf{1}_{\rho \leq \tau}]\leq E^{P}[M_{\rho}] \leq m,
\eeno
which verifies the claim. We then derive that
\beno
E^{P}[V(\tau,\bar{\omega}\otimes_t\omega,M_{\tau})\wedge\epsilon^{-1}]
 &\leq& E^{P}[E^{Q_{\epsilon}^{\omega}}[\xi]]+\epsilon
 = E^{P\otimes_{\tau}Q_{\epsilon}^{\cdot}}[\xi]+\epsilon
 \\&\leq& \sup_{P^{\prime}\in\mathcal{P}(t, \bar{\omega}, m)}E^{P^{\prime}}[\xi]+\epsilon= V(t,\bar{\omega},m) + \epsilon.
\eeno
Let $\epsilon \rightarrow 0$, we have
$E^{P}[V(\tau,\bar{\omega}\otimes_t\omega,M_{\tau})]\leq V(t,\bar{\omega},m).$
As $M\in\mathcal{M}^{+}_{t,\bar{\omega},m}(P)$ is arbitrary, we get
$ \sup_{M\in\mathcal{M}^{+}_{t,\bar{\omega},m}(P)}E^{P}[V(\tau,\bar{\omega}\otimes\omega,M_{\tau})]\leq V(t,\bar{\omega},m).$
In addition, as $P\in\mathcal{P}(t,\bar{\omega},m)$ is arbitrary, we arrive at (\ref{eq:2.5}) by taking supremum over $\mathcal{P}(t,\bar{\omega},m)$, which completes the proof.
\qed

\section{Stochastic Control under Expectation Constraints}\label{sec-cont}
Define $\Omega^{\prime} =  \{\omega^{\prime} \in C([0,T];\mathbb{R}^{n}):\omega^{\prime}_{0} = 0\}$ for some $n\in\mathbb{N}\setminus\{0\}$. Similar to $(\Omega, \mathcal{F}, P_{0})$ in Section \ref{sec-main}, we consider the probability space $(\Omega^{\prime}, \mathcal{F}^{\prime}, P_{0}^{\prime})$ to model real world scenarios. 

\subsection{\textbf{Strong DPP for Stochastic Control Problems}}
For each $(t,\omega)\in [0,T]\times\Omega$, we are given a non-empty set $\mathcal{U}(t,\omega)$ whose elements are interpreted as controls starting from time $t$ with the past path $\omega$. Note that elements in $\Omega$ are observable, thus the dependence of control on past paths is reasonable. We assume that $\mathcal{U}(t,\omega)$ depends on $\omega$ only up to time $t$ in the sense that we do not distinguish two paths at time $t$ if they coincide up to time $t$. That is,
$\mathcal{U}(t,\omega) = \mathcal{U}(t,\tilde{\omega}) ~\text{if}~\omega=\tilde{\omega}~\text{on}~[0,t].$ For each $(t,\omega)\in [0,T]\times\Omega$ and $\nu \in \mathcal{U}(t,\omega) $, we are given a continuous process $X^{t,\omega,\nu} : [0, T]\times \Omega^{\prime} \rightarrow \mathbb{R}^{d}$ satisfying $X^{t,\omega,\nu}_{s}(\omega^{\prime}) = \omega_{s}$ for all $s\in[0,t]$ and $\omega^{\prime} \in \Omega^{\prime}$, which indicates that we can not change the past.

\begin{Assumption}\label{assumfg}
For functions $f(\cdot),g(s,\cdot): \Omega \rightarrow \bar{\mathbb{R}}$, $s\in[0,T]$, let us assume that $f(\cdot)$ is upper semi-analytic and $g(s, \cdot)$ is lower semi-analytic and lower semi-continuous for all $s\in[0,T]$.
\end{Assumption}

Here $g(s,\omega):=g(s,\omega|_{[0,s]})$. For each $(t,\omega,m) \in [0,T]\times\Omega\times\mathbb{R}$, the set of admissible controls with the constraint level $m$ is defined as
\begin{equation} \label{eq:3.1}
\mathcal{U}(t,\omega,m) := \{\nu \in \mathcal{U}(t,\omega): E[g(s,X^{t,\omega,\nu})]\leq m\},\ \ \ \forall s\in[t,T],
\end{equation}
where we denote $E[g(s,X^{t,\omega,\nu})]:=E[g(s, X^{t,\omega,\nu}|_{[0,s]})]$.

We consider the value function
\begin{equation} \label{eq:3.2}
V(t,\omega,m)= \sup_{\nu \in \mathcal{U}(t,\omega,m)} E[f(X^{t,\omega,\nu})].
\end{equation}
Note that expectations in (\ref{eq:3.1}) and (\ref{eq:3.2}) are taken under $P_{0}^{\prime}$.

Let us consider the setting of stochastic control problems with the admissible set $\mathcal{U}(t,\omega)$ such that for each $\nu \in \mathcal{U}(t,\omega)$, the controlled process $X^{t,\omega,\nu}: (\Omega^{\prime}, \mathcal{F}^{\prime}) \rightarrow (\Omega,\mathcal{F})$ is measurable. That is, for every $A \in \mathcal{F}$, we have $(X^{t,\omega,\nu})^{-1}(A) \in \mathcal{F}^{\prime}$.
Note that the process $X^{t,\omega,\nu}$ induces a probability measure $P_{t,\omega,\nu}$ on $(\Omega,\mathcal{F})$ by
$P_{t,\omega,\nu} (A) := P_{0}^{\prime}((X^{t,\omega,\nu})^{-1}(A)),~~A\in\mathcal{F}.$

We call $P_{t,\omega,\nu}$ the probability induced by $\nu$ for any $\nu \in \mathcal{U}(t,\omega)$, and define $\mathcal{P}(t,\omega) = \{P_{t,\omega,\nu}\in\mathfrak{P}(\Omega): \nu \in \mathcal{U}(t,\omega)\}$ as the set of probability measures induced by elements in $\mathcal{U}(t,\omega)$. If $\omega=\tilde{\omega}~\text{on}~[0,t],$ then
$\mathcal{P}(t,\omega) = \mathcal{P}(t,\tilde{\omega}).$
In addition, we have $\mathcal{P}(t,\omega) \neq \emptyset$ for any $(t,\omega) \in [0,T]\times\Omega$, and $P(\Omega^{\omega}_{t}) =1$ for any $P \in \mathcal{P}(t,\omega)$.
By the definition of $P_{t,\omega,\nu}$, let us denote
$$E^{P_{t,\omega,\nu}}[f] := E[f(X^{t,\omega,\nu})]~\text{and}~E^{P_{t,\omega,\nu}}[g(s)] := E[g(s,X^{t,\omega,\nu}|_{[0,s]})].$$
Finally, the set $\mathcal{P}(t,\omega,m)$ of admissible probability measures is defined by
$\mathcal{P}(t,\omega,m) := \{P\in\mathcal{P}(t,\omega):E^{P}[g(s)]\leq m,\ \forall s\in[t,T]\}$. We then have
\begin{equation} \label{eq:3.3}
V(t,\omega,m)  = \sup_{\nu\in\mathcal{U}(t,\omega,m)}E[f(X^{t,\omega,\nu})]
  = \sup_{P\in\mathcal{P}(t,\omega,m)}E^{P}[f].
\end{equation}

Similar to Definition \ref{defM} in the general framework, we need the assistance of some auxiliary supermartingales to establish the DPP result.
\begin{definition}\label{defMM}
For each $\nu \in \mathcal{U}(t,\omega, m)$, let $\mathcal{M}^{+}_{t,\omega,m}(\nu)$ be the collection of all supermartingales on $[t, T] \times \Omega \rightarrow \bar{\mathbb{R}}$ such that\\
$\mathrm{(i)}$ $M_{t} \leq m$ for $P_{0}^{\prime}$-a.e.;\\
$\mathrm{(ii)}$ $M(X^{t,\omega,\nu})$ is a supermartingale under $P_{0}^{\prime}$;\\
$\mathrm{(iii)}$ $M_{s}(X^{t,\omega,\nu}) \geq g(s,X^{t,\omega,\nu})$ for $P_{0}^{\prime}$-a.e..
\end{definition}

The goal is to prove the following DPP
\begin{equation} \label{DPP-for-CPEC}
\begin{split}
V(t,\omega,m) & = \sup_{\nu\in\mathcal{U}(t,\omega,m)}\sup_{M\in\mathcal{M}^{+}_{t,\omega,m}(\nu)}E[V(\tau,X^{t,\omega,\nu},M_{\tau})] \\
 & = \sup_{\nu\in\mathcal{U}(t,\omega,m)}\inf_{M\in\mathcal{M}^{+}_{t,\omega,m}(\nu)}E[V(\tau,X^{t,\omega,\nu},M_{\tau})].
\end{split}
\end{equation}

\begin{theorem} \label{thm:2}
Suppose that $f$ and $g$ satisfy Assumption \ref{assumfg},
 and sets $\mathcal{P}(t,\omega)$ induced by sets $\mathcal{U}(t,\omega)$ satisfy Assumption \ref{ass:0} and Assumption \ref{ass:1}, we have that \eqref{DPP-for-CPEC} holds.
\end{theorem}
{\it Proof}
It is easy to check $\mathcal{P}(t,\omega,m)=\{P_{t,\omega,\nu}\in\mathfrak{P}(\Omega): \nu \in \mathcal{U}(t,\omega,m)\}$ and $\mathcal{M}^{+}_{t,\omega,m}(\nu) = \mathcal{M}^{+}_{t,\omega,m}(P_{t,\omega,\nu})$, so (\ref{DPP-for-CPEC}) follows directly from Theorem \ref{thm:main}.
\qed

\subsection{\textbf{Connection to Classical DPP}}
The DPP in (\ref{DPP-for-CPEC}) involves both the supremum and infimum  over $\mathcal{M}^{+}_{t,\omega,m}(\nu)$, because the set $\mathcal{M}^{+}_{t,\omega,m}(\nu)$ of supermartingales depends on the control $\nu$. In some specific cases such as the next proposition, we can actually get rid of the supremum and infimum over $\mathcal{M}^{+}_{t,\omega,m}(\nu)$.
\begin{proposition} \label{pro:1}
Define 
\begin{align}\label{Mdefin}
\mathcal{M}^{+}_{t,\omega,m} := \cap_{\nu\in\mathcal{U}(t,\omega,m)}\mathcal{M}^{+}_{t,\omega,m}(\nu),
\end{align}
and assume $\mathcal{M}^{+}_{t,\omega,m}$ is non-empty. Under the assumptions in Theorem \ref{thm:2}, for any $M \in \mathcal{M}^{+}_{t,\omega,m}$, it holds that
\begin{equation} \label{Simple-DPP-for-CPEC}
V(t,\omega,m)  = \sup_{\nu\in\mathcal{U}(t,\omega,m)}E[V(\tau,X^{t,\omega,\nu},M_{\tau})].
\end{equation}	
\end{proposition}

Note that \eqref{Simple-DPP-for-CPEC} is a strong version of DPP in a non-Markovian setting. In a Markovian setting, 
one has $V(t,\omega,m) = V(t,\omega_{t},m)$ using $x$ to denote the state variable and
\eqref{Simple-DPP-for-CPEC} is written as
$V(t,x,m)  = \sup_{\nu\in\mathcal{U}(t,x,m)}E[V(\tau,X^{t,x,\nu}_{\tau},M_{\tau})]$,
which is a strong DPP in \cite{Bouchard-Nutz} with generalized state constraints.
Therefore, in the Markovian setting, \eqref{Simple-DPP-for-CPEC}
can be used to show that the value function of a control
problem with general state constraints is the viscosity solution of some associated HJB equation.

\subsection{\textbf{General Dynamic State Constraints}}
We now give a general framework to describe dynamic state constraints. Fix analytic sets $\mathcal{O}(s, \omega) \subseteq \Omega$ indexed by time $(s, \omega)\in [0,T] \times \Omega$. Note that $\mathcal{O}(s, \omega)$ depends only on $s$ and the path $\omega|_{[0,s]}$ up to time $s$.
We are interested in dynamic constraints that the controlled state process $X^{t,\omega,\nu} \in \mathcal{O}(s, X^{t,\omega,\nu})$ at each time $s\in[t, T]$.
More precisely, the admissible control set is defined by
\ben \label{general-state-time-constraint-control}
\bar{\mathcal{U}}(t,\omega) = \{\nu\in\mathcal{U}(t,\omega):  X^{t,\omega,\nu} \in \mathcal{O}(s, X^{t,\omega,\nu}), \ P_{0}^{\prime}\text{-a.s.},\  \forall \ s\in[t,T]\}.
\een
We then study the control problem under the state constraint (CPSC):
\ben \label{control-problem-with-dsc}
\bar{V}(t,\omega)= \sup_{\nu \in \bar{\mathcal{U}}(t,\omega)} E[f(X^{t,\omega,\nu})].
\een
Let us consider the set of all paths satisfying state constraints up to time $t$ by
\ben \label{up-to-time-t}
\Omega(\mathcal{O}(t)) := \{\omega \in \Omega:
\omega \in \mathcal{O}(s, \omega),
 \forall \ s \in [0,t]\},\een
We assume 
$
\bar{\mathcal{U}}(t,\omega)\neq\emptyset,
$
for any $t \in [0,T], \omega \in \Omega(\mathcal{O}(t))$ to exclude the trivial case. In order to rewrite the above dynamic state constraints as expectation constraints, we set the function $g(s,\omega)$ by
\begin{equation} \label{state-g-function-new}
g(s,\omega):=\begin{cases}
               0,&\omega\in \Omega(\mathcal{O}(s)),\\
               1,& \text{otherwise.}
            \end{cases}
\end{equation}
Based on the pair of $f(\cdot)$ and $g(s,\cdot)$, $s\in[0,T]$, we define $\mathcal{U}(t,\omega,m)$ as in \eqref{eq:3.1}  and
$V(t,\omega,m)$ as in \eqref{eq:3.2}.
By the definition of $g$ in \eqref{state-g-function-new}, we have
\begin{equation}
X^{t,\omega,\nu} \in \mathcal{O}(s, X), ~P_{0}^{\prime}\text{-a.s.,} ~\forall s\in[t,T] ~~\Longleftrightarrow ~~E[g(s, X^{t,\omega,\nu})]\leq 0,\ \forall s\in[t,T]. \nonumber
\end{equation}
It follows that
\begin{equation} \label{state-equivalent-2-new}
\bar{\mathcal{U}}(t,\omega) = \mathcal{U}(t,\omega,0) ~\text{and}~ \bar{V}(t,\omega) = V(t,\omega,0).
\end{equation}
For each $\nu \in \mathcal{U}(t,\omega,m)$, $\mathcal{M}^{+}_{t,\omega,m}(\nu)$ and $\mathcal{M}^{+}_{t,\omega,m}$ are defined by the modification of Definition \ref{defMM} and the set in \eqref{Mdefin} with the specific function $g$ in \eqref{state-g-function-new}.

\begin{theorem} \label{general-state-time-constraint}
Suppose the sets $\mathcal{P}(t,\omega)$ induced by $\mathcal{U}(t,\omega)$ satisfy Assumption \ref{ass:1}, and the function $g$ has lower semi-continuous path,
then the corresponding DPP for \eqref{control-problem-with-dsc} holds true that
\begin{equation} \label{DPP-for-general-state-time-constraint}
\bar{V}(t,\omega)=\sup_{\nu \in \bar{\mathcal{U}}(t,\omega)}E[\bar{V}(\tau,X^{t,\omega,\nu})].
\end{equation}
\end{theorem}
{\it Proof.}
By Theorem \ref{thm:2}, the DPP $(\ref{DPP-for-CPEC})$ for $V$ holds true. Take $m =0$, if we have $\nu \in \mathcal{U}(t,\omega,0)$, then for $P_{0}^{\prime}$-a.e., $g(s, X^{t,\omega,\nu}) = 0$ for $s\in[t,T]$. Therefore the constant process $0 \in \mathcal{M}^{+}_{t,\omega,m}$, and 
$
V(t,\omega,0) = \sup_{\nu\in \mathcal{U}(t,\omega,0)}E[V(\tau,X^{t,\omega,\nu},0)] 
$ follows from Proposition \ref{pro:1},
which implies \eqref{DPP-for-general-state-time-constraint} thanks to \eqref{state-equivalent-2-new}.
\qed

\section{Applications in Quantitative Finance}\label{sec-app}
This section is devoted to applications in the context of optimal investment and option hedging under dynamic trading constraints. We aim at reformulating each of those constraints to expectation constraints and consequently DPP for these problems can be established by applying Theorem \ref{general-state-time-constraint}.

%Moreover, in certain models, the value function can be shown as a constrained viscosity solution of some path-dependent PDE as presented earlier.
%To make the presentation easy to follow, we group the first three constraints (state, floor, drawdown) in one subsection as they share similar structure. We leave the other two (target problem, quantile hedging) in a separate subsection target problem can be seen as a special case of quantile hedging.

\subsection{\textbf{State, Floor, Drawdown Constraints}}
\textbf{Case 1}: State Constraint.\\
Let $d\geq 1$ and $\Omega =  \{\omega \in C([0,T];\mathbb{R}^{d}):\omega_{0} = 0\}$. Fix a family of analytic sets $O(t) \subseteq \mathbb{R}^{d}$ indexed by time $t\in[0,T]$, the state constraint condition requires the controlled process to stay in some open sets ($X(t) \in O(t)$) at each intermediate time $t$.
We define 
\begin{align*}
\mathcal{O}_{1}(s, \omega):= \{\omega \in \Omega: \omega_{t} \in O(t), \forall\ t \in [0,s]\}.
\end{align*}
Replacing $\mathcal{O}(s, \omega)$ by $\mathcal{O}_{1}(s, \omega)$ in \eqref{general-state-time-constraint-control} and further replacing $\bar{\mathcal{U}}(t,\omega)$ by $\bar{\mathcal{U}}_{1}(t,\omega)$ in \eqref{control-problem-with-dsc}, we can define $\bar{V}_{1}(t,\omega)$. We refer the resulting problem as control problem under the state constraint (CPSC).
\ \\
\ \\
\textbf{Case 2}: Floor Constraint.\\
Let $d=1$ in this case and $\Omega =  \{\omega \in C([0,T];\mathbb{R}):\omega_{0} = 0\}$. A fixed continuous path $\beta \in \Omega$
is regarded as the floor and it is required that the controlled process stays above $\beta$.
Let $\mathcal{O}_{2}(s, \omega)=\mathcal{O}(s) := \{\omega \in \Omega: ~\omega_{t} \geq \beta_t, \forall \ t \in [0,s]\}$.
Replacing $\mathcal{O}(s, \omega)$ by $\mathcal{O}_{2}(s, \omega)$, we
have the admissible control set $\bar{\mathcal{U}}_{2}(t,\omega)$ by \eqref{general-state-time-constraint-control}
and value function $\bar{V}_{2}(t,\omega)$ by \eqref{control-problem-with-dsc}.
We refer the resulting problem as control problem
under floor constraint (CPFC).\\
\ \\
\textbf{Case 3}: Drawdown Constraint.\\
Let $d=1$ in this case. Fix $x\geq 0$, and let $\Omega =  \{\omega \in C([0,T];\mathbb{R}):\omega_{0} = x\}$. Fix a continuous function $\alpha:
[0,T] \rightarrow [0,1]$, the controlled state process is required to satisfy the drawdown condition at all intermediate time (with drawdown no less than $1-\alpha$). Given $\omega \in \Omega$, we define its running maximum function by
$\omega^{*}_{s} = \sup_{0 \leq r \leq s} \omega_{r}.$
Let us define $\mathcal{O}_{3}(s, \omega) := \{\omega \in \Omega: \omega_{t} \geq \alpha(t)\omega^{*}_{t}, \forall
\ t \in [0,s]\}$.
Replacing $\mathcal{O}(s, \omega)$ by $\mathcal{O}_{3}(s, \omega)$, we
have the admissible control set $\bar{\mathcal{U}}_{3}(t,\omega)$ by \eqref{general-state-time-constraint-control}
and value function $\bar{V}_{3}(t,\omega)$ by \eqref{control-problem-with-dsc}.
We refer the resulting problem as control problem
under drawdown constraint (CPDC).

Replacing $\mathcal{O}(s, \omega)$ by $\mathcal{O}_{i}(s, \omega)$, we define $\Omega(\mathcal{O}_{i}(t))$ by \eqref{up-to-time-t} and $g_{i}$ by \eqref{state-g-function-new} for  $i = 1,2,3$. All three functions $g_{i}~(i = 1,2,3)$ are lower semi-analytic because the associated  $O(t)$ is analytic, the path $\beta$ is continuous and the function $\alpha$ is continuous. The next result is a direct application of Theorem \ref{general-state-time-constraint}.
\begin{theorem} \label{thm:3}
Suppose that the sets $\mathcal{P}(t,\omega)$ induced by $\mathcal{U}(t,\omega)$ satisfy Assumption \ref{ass:1}, and the functions $g_{i}~(i = 1,2,3)$ have lower semi-continuous path,
then DPP results for \textup{(CPSC), (CPFC)} and \textup{(CPDC)} hold that
\begin{equation} \label{DPP-for-s-f-d}
\bar{V}_{i}(t,\omega)=\sup_{\nu \in \bar{\mathcal{U}}_{i}(t,\omega)}E[\bar{V}_{i}(\tau,X^{t,\omega,\nu})],\ \ \text{for}\ i = 1,2,3.
\end{equation}
\end{theorem}

\begin{remark} In case 1 and case 2, it is by no means restrictive to assume that the starting point of the controlled process is at the origin. If not, we can simply make the translation such that it starts from origin. However, case 3 is different. If a path satisfies the maximum drawdown condition, it may not satisfy the condition after a translation.  Therefore, in case 3, we allow the controlled process to start from any point $x \geq 0$.
\end{remark}

\subsection{\textbf{Target Problem and Quantile Hedging}}
We consider in this subsection the quantile hedging and its implication to target constraints. Fix a quantile level
$\gamma \in [0,1]$ and a family of analytic target sets $G(t) \subseteq \mathbb{R}^{d}$ indexed by time $t\in[0,T]$. The quantile hedging problem requires the probability to be greater than $\gamma$ that the controlled process stays in the target sets ($X(t) \in G(t)$). We consider the control problem under quantile hedging constraints (CPQC) that 
\begin{align}\label{defbarV4}
\bar{V}_{4}(t,\omega, m)= \sup_{\nu \in \bar{\mathcal{U}}_{4}(t,\omega, m)} E[f(X^{t,\omega,\nu})]
\end{align}
where 
\begin{align*}
\bar{\mathcal{U}}_{4}(t,\omega, m):= \{\nu\in\mathcal{U}(t,\omega):P_{0}^{\prime}(X_{s}^{t,\omega,\nu} \in G(s))\geq m,\ s\in[0,T]\}.
\end{align*}
We can transform dynamic quantile hedging constraints to expectation constraints on paths by setting $g_{4}$ as
\begin{equation}
  g_{4}(s, \omega)=\begin{cases}
               0,&\omega_{s}\in G(s),\\
               1,& \text{otherwise.}\nonumber
            \end{cases}
\end{equation}
Note that $g_{4}$ is lower semi-analytic as  $G(s)$ is analytic.
Based on this pair of $f$ and $g_{4}$, we can define $\mathcal{U}_{4}(t,\omega,m)$ as in \eqref{eq:3.1} and
$V_{4}(t,\omega,m)$ as in \eqref{eq:3.2}.
The definition of $g_{4}$ implies that
\begin{equation}
P_{0}^{\prime}(X_{s}^{t,\omega,\nu} \in G(s))\geq m,\ \ s\in[t,T] ~~\Longleftrightarrow ~~E[g_{4}(s, X^{t,\omega,\nu})]\leq 1-m,\ \ s\in[t,T].\nonumber
\end{equation}
It then follows that
\begin{equation} \label{quantile-equivalent-2}
\bar{\mathcal{U}}_{4}(t,\omega,m) = \mathcal{U}_{4}(t,\omega,1-m) ~\text{and}~ \bar{V}_{4}(t,\omega, m) = V_{4}(t,\omega,1-m).
\end{equation}

\begin{definition}
For each $\nu \in \bar{\mathcal{U}}_{4}(t,\omega, m)$, let $\bar{\mathcal{M}}_{t,\omega,m}^{-}(v)$ be the collection of all submartingales on $[t, T] \times \Omega \rightarrow \bar{\mathbb{R}}$ such that\\
$\mathrm{(i)}$ $M_{t} \geq m$ for $P_{0}^{\prime}$-a.e.;\\
$\mathrm{(ii)}$ $M(X^{t,\omega,\nu})$ is a submartingale under $P_{0}^{\prime}$;\\
$\mathrm{(iii)}$ $M_{s}(X^{t,\omega,\nu}) \leq 1 - g_{4}(s, X^{t,\omega,\nu})$ for $P_{0}^{\prime}$-a.e. and $s\in[t,T]$.
\end{definition}

Let $\mathcal{M}_{t,\omega,1-m,4}^{+}(v)$ be defined as in Definition \ref{defMM} by using the function $g_{4}$ for $\nu \in \mathcal{U}_{4}(t,\omega,1-m)$, then we have
\ben \label{relation-plus-minus}
\bar{\mathcal{M}}_{t,\omega,m}^{-}(v) = 1 - \mathcal{M}_{t,\omega,1-m,4}^{+}(v) := \{1-M: M \in \mathcal{M}_{t,\omega,1-m,4}^{+}(v)\}.
\een

The next result gives DPP for (CPQC).
\begin{theorem} \label{quantile}
Suppose the sets $\mathcal{P}(t,\omega)$ induced by $\mathcal{U}(t,\omega)$ satisfy Assumption \ref{ass:1}, and the function $g_{4}$ has lower semi-continuous path, then DPP for \textup{(CPQC)} holds that
\begin{equation} \label{DPP-for-qh}
\begin{split}
\bar{V}_{4}(t,\omega, m) & = \sup_{\nu\in\bar{\mathcal{U}}_{4}(t,\omega, m)}\sup_{M\in\bar{\mathcal{M}}_{t,\omega,m}^{-}(v)}E[\bar{V}_{4}(\tau,X^{t,\omega,\nu},M_{\tau})] \\
 & = \sup_{\nu\in\bar{\mathcal{U}}_{4}(t,\omega, m)}\inf_{M\in\bar{\mathcal{M}}_{t,\omega,m}^{-}(v)}E[\bar{V}_{4}(\tau,X^{t,\omega,\nu},M_{\tau})].
\end{split}
\end{equation}
\end{theorem}
{\it Proof.}
Thanks to Theorem \ref{thm:2}, DPP in (\ref{DPP-for-CPEC}) for $V_{4}$ holds valid. That is, we have
\begin{equation} \label{DPP-for-CPEC-TPQC}
\begin{split}
V_{4}(t,\omega,1-m) & = \sup_{\nu\in\mathcal{U}_{4}(t,\omega, 1-m)}\sup_{M\in\mathcal{M}_{t,\omega,1-m,4}^{+}(v)}E[V_{4}(\tau,X^{t,\omega,\nu},M_{\tau})] \\
 & = \sup_{\nu\in\mathcal{U}_{4}(t,\omega,1-m)}\inf_{M\in\mathcal{M}_{t,\omega,1-m,4}^{+}(v)}E[V_{4}(\tau,X^{t,\omega,\nu},M_{\tau})].
\end{split}
\end{equation}
Together with \eqref{quantile-equivalent-2} and \eqref{relation-plus-minus}, we can deduce (\ref{DPP-for-qh}).
\qed

By choosing $m=1$, it is clear that the constant process $1 \in \cap_{\nu\in\bar{\mathcal{U}}_{4}(t,\omega, 1)} \bar{\mathcal{M}}_{t,\omega,1}^{-}(v)$, and we obtain DPP for target constraint in the next result.
\begin{proposition} \label{target-as-a-proposition}
Under assumptions in Theorem \ref{quantile}, we have
\begin{equation} \label{Simple-DPP-for-CPTC}
\bar{V}_{4}(t,\omega, 1)  = \sup_{\nu\in\bar{\mathcal{U}}_{4}(t,\omega, 1)}E[\bar{V}_{4}(\tau,X^{t,\omega,\nu},1)].
\end{equation}	
\end{proposition}

\begin{remark} Proposition \ref{target-as-a-proposition} reduces to a geometric type of DPP when $f$ is chosen to be a constant. In particular, if we define the reachability set:
$$D(t) := \{\omega\in\Omega: \exists ~ \nu \in \mathcal{U}(t,\omega),~\text{such that}~ X^{t,\omega,\nu}_{s}\in G(s) ~P_{0}^{\prime}\text{-a.s.}, \ s\in[t,T]\},$$
the geometric type of DPP states that
\begin{equation} \label{Geometric-DPP}
D(t) = \{\omega\in\Omega: \exists ~ \nu \in \mathcal{U}(t,\omega),~\text{such that}~ X_{\tau}^{t,\omega,\nu}\in D(\tau) ~P_{0}^{\prime}\text{-a.s.}\}.
\end{equation}
By setting $f(x)\equiv 1$ in the problem \eqref{defbarV4}, \eqref{Geometric-DPP} follows directly from \eqref{Simple-DPP-for-CPTC}.
\end{remark}

\section{Conclusions}\label{section_conc}
\begin{onehalfspacing}
We investigated DPP for stochastic control problems under intermediate expectation constraints at each time in a general non-Markovian framework. For continuous state processes, DPP is established in the present paper by applying the measurable selection method together with some auxiliary aggregated supermartingales. Moreover, we show that several types of dynamic trading constraints from financial applications can be transformed into the expectation constraints at each intermediate time $t$. Hence DPP holds in optimal investment and hedging problems under these trading constraints. It is still an open problem whether the value function of the control problem corresponds to a constrained viscosity solution to some path-dependent partial differential equation in a non-Markovian setting, which is left for future research.
\end{onehalfspacing}

\begin{singlespacing}
\begin{acknowledgements}
This work is partially done when the first author was a PhD student in Hong Kong Baptist University.
Yuk-Loong Chow is supported by the Fundamental Research Funds for the Central Universities under the grant 19lgpy242.
Xiang Yu is supported by the Hong Kong Early Career Scheme under grant no. 25302116. Chao Zhou is supported by Singapore MOE (Ministry of Education's) AcRF grant R-146-000-219-112.
\end{acknowledgements}
\end{singlespacing}

\bibliographystyle{spmpsci_unsrt}

\end{document}